\documentclass[a4paper,10pt]{article}

\usepackage{amsthm,amssymb,amsfonts,amscd,amsmath,mathtools,mathrsfs}
\usepackage{stmaryrd}
\usepackage{csquotes}
\usepackage{hyperref}
\usepackage{tikz}
\usetikzlibrary{decorations.pathreplacing,matrix,calc,arrows,shapes.geometric}
\usepackage{listings}

\usepackage[backend=biber,isbn=false,giveninits=true,url=false,eprint=false,doi=false,sortcites]{biblatex}
\newcommand{\dlgraffa}{\{ \hspace{-0.1cm} \{}
\newcommand{\drgraffa}{\} \hspace{-0.1cm} \}}
\newcommand{\vertiii}[1]{{\left\vert\kern-0.25ex\left\vert\kern-0.25ex\left\vert #1 \right\vert\kern-0.25ex\right\vert\kern-0.25ex\right\vert}}

\newcommand{\IR}{{\mathbb R}}
\newcommand{\IC}{{\mathbb C}}

\newcommand{\IZ}{{\mathbb Z}}
\newcommand{\bx}{{\mathbf x}}
\newcommand{\bd}{{\mathbf d}}
\newcommand{\bn}{{\mathbf n}}

\newcommand{\bT}{{\mathbf T}}
\newcommand{\Mh}{{\mathcal M_h}}
\newcommand{\hankel}[1]{{H_{#1}^{(1)}}}

\newcommand{\ri}{{\mathrm i}}
\newcommand{\re}{{\mathrm e}}
\newcommand{\rd}{{\,\mathrm d}}
\newcommand{\inc}{{\mathrm{inc}}}
\newcommand{\tot}{{\mathrm{tot}}}
\newcommand{\scat}{{\mathrm{scat}}}

\newcommand{\pa}{{\mathtt{a}}}
\newcommand{\pb}{{\mathtt{b}}}
\newcommand{\pd}{{\mathtt{d}}}
\newcommand{\dn}{{\partial_\bn}}

\newcommand{\GR}{{\Gamma_R}}
\newcommand{\OmOut}{{\Omega_{\text{o}}}}
\newcommand{\uOut}{{u_{\text{o}}}}
\newcommand{\OmIn}{{\Omega_{\text{i}}}}
\newcommand{\uIn}{{u_{\text{i}}}}
\newcommand{\nI}{{n_{\text{i}}}}
\newcommand{\nO}{{n_{\text{o}}}}
\newcommand{\kO}{{\kappa_{\text{o}}}}
\newcommand{\kI}{{\kappa_{\text{i}}}}

\usepackage{soul} 
\usepackage{xcolor}

\graphicspath{{figs/}}

\title{TMATDG: applying TDG methods to multiple scattering via T-matrix approximation} 
\author{Armando Maria Monforte\thanks{Department of Mathematics, University of Pavia, Italy
		(\href{mailto:armandomaria.monforte01@universitadipavia.it}{armandomaria.monforte01@universitadipavia.it}), ORCID: 0009-0000-7687-2217}
} 
\date{\today} 

\begin{document}
	\newtheorem{thm}{Theorem}[section]
	\newtheorem{prop}[thm]{Proposition}
	\newtheorem{lem}[thm]{Lemma}
	\newtheorem{defin}[thm]{Definition}
	\newtheorem{rem}[thm]{Remark}
	\allowdisplaybreaks
	
	\maketitle
	
	\section*{Abstract}
	We present a \textsc{matlab} package for the solution of multiple scattering problems, coupling Trefftz {Discontinuous} Galerkin methods for Helmholtz scattering with the T-matrix method. The method applies to both impenetrable and penetrable scatterers and relies on plane wave discrete spaces; the DtN-TDG formulation for the transmission problem is new. We rely on the \textsc{tmatrom} package to numerically approximate the T-matrices and deal with multiple scattering problem, providing a framework to handle scattering by polygonal obstacles.
	
	\bigskip\noindent
	{\bf Keywords:} \; Multiple scattering, Helmholtz equation, T-matrix, Discontinuous Galerkin, Trefftz method, Plane wave basis
	
	\bigskip\noindent
	{\bf Mathematics Subject Classification (2020):} \;
	65N30,
	35J05,  
	78A45,
	65N35
	
	\section{Introduction}
	We present \textsc{tmatdg}, a \textsc{matlab} package for the numerical solution of two-dimensional multiple scattering problems of time-harmonic acoustic waves, governed by the Helmholtz equation. The approach is based on the T-matrix framework \cite{Waterman1965}, which reduces a multiple scattering problem to the solution of independent single scattering problems, whose results are then coupled through linear relations. A key feature of the approach is that the T-matrix, which represents the scattering properties of a single obstacle, depends only on the obstacle shape and material parameters and can therefore be computed once, stored, and reused, yielding substantial computational savings for ensembles containing repeated shapes.
	
	We propose the use of a Dirichlet-to-Neumann Trefftz Discontinuous Galerkin (DtN-TDG) method \cite{Kapita2018} to numerically approximate the T-matrix of a single polygonal obstacle. The DtN-TDG method is used to compute the far-field responses generated by a basis of regular incident fields, from which the T-matrix is constructed. The method applies to both impenetrable and penetrable scatterers within the same unified framework and relies on plane wave discrete spaces; while the formulation for impenetrable obstacles is found in literature, the DtN-TDG formulation for the transmission problem is new.
	
	\textsc{tmatdg} is based on the \textsc{tmatrom} package \cite{Ganesh2017} and provides tools for computing, storing, and reusing T-matrices, as well as for solving and visualizing multiple scattering problems, offering an efficient and flexible platform for multiple scattering simulations.
	
	\section{Multiple scattering problem and T-matrix}
	Multiple scattering is the repeated deflection of waves as they interact with many obstacles in a medium \cite{Martin2006}. The field scattered from one obstacle will induce further scattered fields from all the other obstacles, and so on, leading to complex phenomena. 
	
	We consider the scattering of time harmonic acoustic waves by an ensemble of two dimensional obstacles $D_j \subset \IR^2$, with $j=1,\ldots,N$. We consider linear {acoustics} with $\re^{-\ri \omega t}$ dependence on time $t$, where $\ri = \sqrt{-1}$ and $\omega$ is the wave angular frequency. An incident wave $u^\inc$ interacts with the obstacles and induces a scattered field $u^\scat$. The incident field satisfies the Helmholtz equation $\Delta u^\inc +\kappa^2 u^\inc = 0$.
	We both consider impenetrable and homogeneous penetrable obstacles in the ensemble, so $u^\scat$ satisfies the Helmholtz equation with piecewise-constant wavenumber $\kappa \in \IC$.
	
	A common strategy to solve a multiple scattering problem is to deal with every obstacle separately, solving many single scattering problems in an iterative way, coupling the scatterers with one another. 
	
	\subsection{The T-matrix method}
	We present the T-matrix method for multiple scattering; this method was first formulated by Waterman \cite{Waterman1965, Waterman1971} {as a method for scattering by a single particle}. The T-matrix method is particularly powerful for multiple scattering problems \cite{Martin2006}, since it allows to easily couple scatterers. Each scatterer interacts both with the incident wave $u^\inc$ and the wave scattered by all the other obstacles. The translation operator \cite{Waterman1971} is used to couple scatterers taking into account their mutual position, effectively representing the field scattered by one obstacle as a field incident on the other one; all these one-to-one interactions are combined into one system of linear relations, which can be solved to get the total scattered field.
	
	Since the T-matrix method allows us to reduce the study to a single scattering problem, we first focus on the approximation of the field scattered by a two-dimensional obstacle $D$, either impenetrable or penetrable by waves, centered in the origin, contained in a disk of radius $R_D>0$. We can assume that $D$ is centered in the origin since the the translation operator easily allows to translate its T-matrix. We denote by $(r, \theta)$ the polar coordinates centered in the origin.
	The T-matrix method is based on the expansion of the incident and scattered field in the so-called regular and radiating wavefunctions $\psi_l$ and $\phi_m$ respectively:
	\begin{align*}
		& \psi_l(r, \theta) = J_{{|l|}}(\kappa r)\re^{\ri l\theta}, \\
		& \phi_m(r, \theta) = \hankel{{|m|}}(\kappa r)\re^{\ri m\theta}, \qquad \text{for}\,\, l,m \in \IZ.
	\end{align*}
	Expanding $u^\inc$ and $u^\scat$ in these basis we get:
	\begin{align} \label{eq:expansion}
		& u^\inc(r,\theta)=\sum_{l \in \IZ} a_l J_{{|l|}}(\kappa r)\re^{\ri l\theta}, \qquad r \geq 0, \\ \nonumber
		& u^\scat(r,\theta)=\sum_{m \in \IZ} b_m \hankel{{|m|}}(\kappa r)\re^{\ri m\theta}, \qquad r \geq R_D.
	\end{align}
	For common incident fields, such as plane or circular waves, the expansion coefficients $a_l$ are known and given explicitly \cite{Ganesh2012}. The unknowns are the expansion coefficients $b_m$ of the scattered field. The relationship between the column vector $ \mathbf{a} $ of incident field expansion coefficients and the column vector $ \mathbf{b} $ of scattered field expansion coefficients is linear and can be expressed through the infinite transition matrix $\mathbf{T}$ as
	\begin{equation*}
		\mathbf{b} = \bT \mathbf{a}.
	\end{equation*} 
	The elements $T_{ml}$ of this matrix, that is called the T-matrix, describe how each component of the incident wave contributes to each component of the scattered wave. The T-matrix entries depend only on the shape of the scatterer and the wavenumber, but not on the incident field \cite{Waterman1965}. This makes the T-matrix highly reusable, since, once computed, it can be used to compute the scattered field for any given incident field, assuming that one is able to derive the regular wavefunctions expansion.
	
	The T-matrix satisfies the following symmetry relation \cite{Ganesh2009}:
	\begin{equation} \label{sym}
		\bT + \bT^* -2\bT \bT^* = 0,
	\end{equation}
	where $\bT^*$ is the conjugate transpose of $\bT$. {This property holds for both penetrable and impenetrable obstacles, but only when the wavenumber $\kappa \in \IR$.}
	
	The T-matrix is by definition infinite-dimensional, but in numerical applications it has to be numerically truncated. A common rule for choosing the truncation order $N_{\text{TR}}$ of the T-matrix for a scatterer of radius $R_D$ and wavenumber $\kappa$ is \cite{Wiscombe1980,Ganesh2012}:
	\begin{equation*}
		N_{\text{TR}} = \kappa R_D + 4(\kappa R_D)^{1/3} + 5.
	\end{equation*} 
	To compute the T-matrix entries there exists various numerical formulations, such as the extended boundary condition method formulated by Waterman \cite{Waterman1969}, which has some numerical limitations \cite{Mishchenko1994,Mishchenko1996}, or the point matching method \cite{Farafonov2010,HellmersJens2011}. In the following we use the far-field based method {proposed} by Ganesh and Hawkins \cite{Ganesh2008,Ganesh2017}, in which the matrix entries are computed as
	\begin{equation} \label{eq:farfield}
		T_{ml} = \frac14 \sqrt{\frac\kappa\pi} \ri^{|m|}(\ri + 1) \int_{\partial B} \phi_l^\infty \re^{-\ri m \theta} \rd s,
	\end{equation}
	where $B$ is the unit circle and $\phi_l^\infty$ is the numerically approximated far field of the scattered wave generated when choosing as incident field the regular wavefunction $\psi_l$.
	In general, when the truncated T-matrix is computed numerically, the symmetry relation \eqref{sym} can be used to estimate the error.
	
	We also point out how the T-matrix behaves under rotations of the scatterer. Let $D'$ be the scatterer obtained {by} rotating $D$ by a certain angle $\alpha>0$. We {show here how to derive the well-known} relation beetwen the T-matrix $\bT'$ of the rotated scatterer and the original T-matrix $\bT$. The rotated polar coordinates are $(r', \theta')$, with $r'=r$, $\theta' = \theta + \alpha$. When we expand $u^\inc$ and $u^\scat$ we get:
	\begin{align*}
		& u^\inc(r,\theta)=\sum_{l \in \IZ} a_l J_{{|l|}}(kr)\re^{\ri l\theta} =\sum_{l \in \IZ} a_l J_{{|l|}}(kr)\re^{\ri l(\theta' - \alpha)} = \sum_{l \in \IZ} a_l J_{{|l|}}(kr)\re^{\ri l\theta'} \re^{-\ri l \alpha}, \\
		& u^\scat(r,\theta)=\sum_{m \in \IZ} b_m \hankel{{|m|}}(kr)\re^{\ri m\theta}= \sum_{m \in \IZ} b_m \hankel{{|m|}}(kr)\re^{\ri m\theta'} \re^{-\ri m \alpha}.
	\end{align*}
	If we express the relation between the expansion coefficients we get that
	\begin{equation*}
		b'_m =  \re^{-\ri m \alpha} b_m = \sum_{l \in \IZ} \bT'_{ml} \re^{-\ri l \alpha} a_l = \sum_{l \in \IZ} \bT'_{ml} a'_l, 
	\end{equation*}
	so we can write $\mathbf{b}' = \bT' \mathbf{a}'$, where
	\begin{equation} \label{eq:rotation}
		\bT' = \mathbf{D}^* \bT \mathbf{D},
	\end{equation}
	and $\mathbf{D} = \text{diag}(\re^{\ri l \alpha})_{l \in \IZ}$ is a diagonal matrix.
	
	\subsection{The TMATROM package}
	When dealing with multiple obstacles, it is computationally convenient to compute the T-matrices individually for each obstacle, possibly in parallel, and then couple them to get the total field. This also allows for easy changes in the ensemble configuration. The \textsc{tmatrom} package \cite{Ganesh2017}, developed by Ganesh and Hawkins, is an object-oriented \textsc{matlab} package for computing T-matrices in two dimensions. It relies on the far-field formulation of the T-matrix \eqref{eq:farfield} and allows to include any type of numerical method to approximate the solution of the single scattering problem, assuming that it is able to compute the far field in a given direction. 
	
	\textsc{tmatrom} allows to easily deal with multiple scattering problems: the built-in functions allow  to put different T-matrices in relation and to solve the system corresponding to the multiple scattering configuration in an iterative way using the GMRES algorithm \cite{Ganesh2013}. The T-matrix expansion is valid as long as the disks containing the {scatterers} do not overlap, otherwise the problem is numerically unstable. This package is capable of dealing with challenging problems such as metamaterial simulations \cite{Hawkins2024}.
	
	Using this package, we are able to include our DtN-TDG method to efficiently approximate the solution of the single scattering problem and so the corresponding T-matrix, both in the penetrable and impenetrable case.

	\section{DtN-TDG approximation of the single scattering problem}	
	We present the DtN-TDG method for the solution of the single scattering problem: this method was first introduced by Kapita and Monk in \cite{Kapita2018} for the Dirichlet obstacle, while the formulation for a penetrable obstacle is new.
	
	\begin{figure}[htbp]
		\centering
		\includegraphics[width=0.25\textwidth]{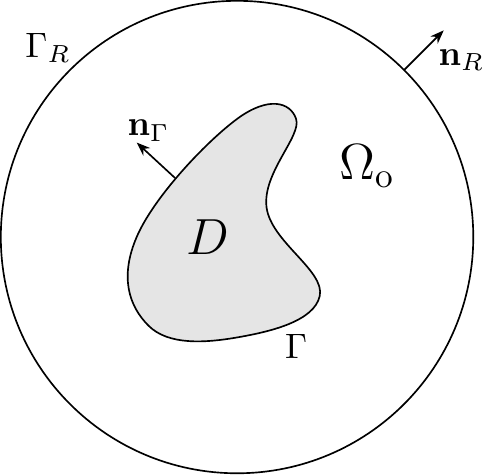}\hfil\includegraphics[width=0.25\textwidth]{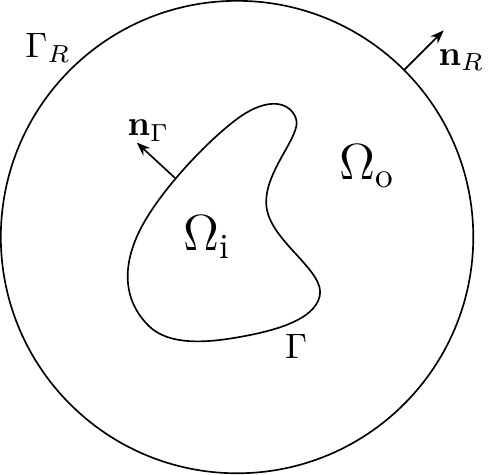}
		\caption{Domain geometry for Dirichlet (left) and transmission (right) problem.}
		\label{fig:domain}
	\end{figure} 
	
	Let $D \subset \IR^2$ be, as before, a Lipschitz obstacle centered {at} the origin and with radius $R_D > 0$, either penetrable or impenetrable by waves. Let $\Gamma := \partial D$ and assume that the refraction index outside of $D$ is constant and equal to $n = \nO$. An incident field $u^\inc$, satisfying the Helmholtz equation $\Delta u^\inc +k^2\nO u^\inc = 0$ interacts with $D$. We consider two instances:
	\begin{itemize}
		\item[(i)] Sound-soft obstacle: we impose a Dirichlet condition $u=-u^\inc$ on $\Gamma$ and seek for $u=u^\scat$ that satisfies the Helmholtz equation in $\IR^2 \setminus D$;
		\item[(ii)] Penetrable obstacle: we assume that $D$ is penetrable by waves with constant refraction index $\nI \in \IC$; we define $u=\uOut=u^\scat$ in $\IR^2 \setminus D$ and $u=\uIn=u^\tot$ in $D$ and impose a transmission condition on $\Gamma$, asking that $u$ solves the Helmholtz equation with variable wavenumber.
	\end{itemize}
	
	To approximate the solution, we introduce an artificial circular boundary $\GR = \partial B_R$ enclosing $D$, such that $R > R_D$. In both problem (i) and (ii), we define the outside domain as $\OmOut := B_R \setminus D$; in problem (i), we solve the Helmholtz equation in $\Omega := \OmOut$, while in problem (ii) we rename $\OmIn:= D$, and define $\Omega := \OmIn \cup \OmOut$. We define $\kO := k \sqrt{\nO}$ and $\kI := k \sqrt{\nI}$, so that $\kappa = k\sqrt{n} \in \IC$ is piecewise constant in $\Omega$; we also denote as $\bn_R$ the outward pointing normal to $\Gamma_R$ and as $\bn_\Gamma$ the normal to $\Gamma$ pointing outside of $D$.
	The problem geometry is depicted in Figure \ref{fig:domain}.
	
	To impose a radiation condition on $\GR$, we make use of a Dirichlet-to-Neumann (DtN) operator $T$, which is defined as
	\begin{align} \label{DtN}
		& T:H^{1/2}(\Gamma_R) \to H^{-1/2}(\Gamma_R),\\ \nonumber
		& Tw(R,\theta) := \kO \sum_{l \in \IZ}w_l\frac{H_l^{(1)'}(\kO R)}{H_l^{(1)}(\kO R)}\re^{\ri l \theta}, \qquad \text{for} \,\, w(R,\theta) = \sum_{l \in \IZ} w_l \re^{\ri l \theta}.
	\end{align}
	Here we denote as $(r,\theta)$ the polar coordinates centered in the origin. The expression of the DtN operator comes from asking that $u^\scat$ satisfies the radiation condition as $r \to \infty$. In particular, if $w$ is a solution of the scattering problem, we have that $Tw = \partial_r w$ on $\GR$. The derivation of $T$ is described in \cite{Kapita2018}, with the only difference being that it makes use of Hankel functions of the second kind, since the radiation condition is taken with the opposite sign.
	
	We use the DtN operator to define the following boundary value problems:
	\begin{align} \label{bvp}
		\begin{cases}
			\Delta u +\kappa_{\text{o}}^2 u = 0 & \text{in} \,\, \Omega, \\
			u = g_D & \text{on} \,\, \Gamma, \\
			\partial_{\bn_R} u - Tu = 0 & \text{on} \,\, \Gamma_R,
		\end{cases}
		& & 
		\begin{cases}
			\Delta u_o +\kappa_{\text{o}}^2 \, u_o = 0 & \text{in} \,\, \OmOut, \\
			\Delta u_i +\kappa_{\text{i}}^2 \, u_i = 0 & \text{in} \,\, \OmIn, \\
			u_i = u_o - g_D & \text{on} \,\, \Gamma, \\
			\partial_{\bn_\Gamma} u_i = \partial_{\bn_\Gamma} u_o + g_N & \text{on} \,\, \Gamma, \\
			\partial_{\bn_R} u_o - Tu_o = 0 & \text{on} \,\, \Gamma_R ,
		\end{cases}
	\end{align}
	where $g_D = -u^\inc$, $g_N = \partial_{\bn_\Gamma} u^\inc$.
	
	Problem \eqref{bvp} is well-posed (see \cite{Kapita2018}), but the definition of the DtN operator \eqref{DtN} includes an infinite series. In practical computations, we truncate the infinite series to a certain $M>0$, thus defining a truncated DtN operator $T_M$. 
	If $M$ is sufficiently large, the truncated BVP with Dirichlet boundary conditions is still well-posed \cite[Thm.~1]{Kapita2018}; similar arguments can be used to prove the well-posedness of the transmission problem.
	
	\subsection{DtN-TDG approximation on the circle}	
	We present the numerical formulation of the DtN-TDG method.	We extend the method described in \cite{Kapita2018} to also solve the transmission problem, adapting the numerical fluxes on $\Gamma$. 
	
	From now on, we assume that $D$ is a polygon and that the relative permittivity assumes a constant value $\nO$ outside of $D$ and possibly a different constant value $\nI$ inside $D$.
	Given a mesh $\mathcal{M}_h = \{ K \}$ of $\Omega$, which has to be such that the wavenumber is constant in every element, we denote by $F_h = \bigcup_{K \in \mathcal{M}_h } \partial K $ the skeleton of the mesh and by $F_h^I = F_h \setminus(\Gamma\cup\GR)$ its inner part. We allow the mesh to have curved edges in correspondence to the circular boundary $\GR$. We then define the main space where our basis functions live, the Trefftz space $T(\mathcal{M}_h)$:
	\begin{align} \label{broken_sob}
		T(\mathcal{M}_h) := \{ & v \in L^2(\Omega) : v_{| K}\in H^1(K) \hspace{0.2cm} \forall K \in \mathcal{M}_h, \\
		& \Delta v + \kappa^2v = 0 \hspace{0.1cm} \text{in} \hspace{0.1cm} K \hspace{0.1cm} \text{and} \hspace{0.1cm} \partial_{\mathbf{n}_K}v \in L^2(\partial K) \hspace{0.2cm}\forall K \in \mathcal{M}_h \}. \nonumber
	\end{align}
	The discrete Trefftz space $V_p(\mathcal{M}_h)$ is a finite-dimensional subspace of $T(\mathcal{M}_h)$ and can be represented as $V_p(\mathcal{M}_h) = \bigoplus_{K \in \mathcal{M}_h}V_{p_K}(K)$, where $V_{p_K}(K)$ is a $p_K$-dimensional subspace of $T(\mathcal{M}_h)$ of functions supported in $K$.
	
	We can then proceed as in \cite{Kapita2018, Hiptmair2011} and derive the TDG formulation. We multiply the Helmholtz equation by a test function $v$ and integrate by parts twice on each $K\in\mathcal{M}_h$. We then replace $u$ and $v$ by discrete functions $u_h, v_h \in V_p(\mathcal{M}_h)$; on $\partial K$ we replace the trace of $u_h$ and $\nabla u_h$ by the numerical fluxes $\hat{u}_h$ and $-\widehat{\ri\kappa\sigma}_h$, that are single-valued approximations of $u_h$ and $\nabla u_h$, respectively, on each edge, obtaining:
	\begin{equation} \label{TDG}
		\int_{\partial K} \hat{u}_h \, \overline{\partial_{\mathbf{n}_K} v_h}  \rd s 
		+ \int_{\partial K} \widehat{\ri\kappa\sigma}_h \cdot \mathbf{n_K} \, \overline{v_h}  \rd s = 0.
	\end{equation}
	
	On internal faces the definition of the fluxes is the classical one \cite{Hiptmair2011}:
	\begin{equation*}
		\begin{cases}
			\hat{u}_h = \dlgraffa u_h\drgraffa  - \pb\ri \kappa^{-1} \llbracket \nabla_h u_h  \rrbracket_N, & \\
			\widehat{\ri\kappa\sigma}_h = - \dlgraffa \nabla_h u_h\drgraffa  - \ri \kappa \pa \llbracket u_h  \rrbracket_N, & 
		\end{cases}
	\end{equation*}
	where  ${\pa, \pb} \in L^\infty(F_h^I \cup \Gamma)$ are positive flux coefficients and $\kappa = \kappa_i$ on $\OmIn$ and $\kappa= \kappa_o$ on $\OmOut$. On the circular boundary $\GR$ we define the fluxes as:
	\begin{equation*}
		\begin{cases}
			\hat{u}_h = u_h - \pd\ri \kappa^{-1} \left(\nabla_h  u_h \cdot \bn_R - T_M u_h \right), & \\
			\widehat{\ri \kappa \sigma}_h = -  T_M u_h \bn_R + \pd\ri \kappa^{-1} T_M^{*} \left( \nabla_h  u_h -T_M u_h \bn_R \right), &
		\end{cases}
	\end{equation*}
	where ${\pd} \in L^\infty(\GR)$ is a positive flux coefficient and $T^{*}$ is the $L^2(\GR)$-adjoint of $T$,
	while on $\Gamma$ we distinguish between the Dirichlet (i) and transmission (ii) problem:
	\begin{align*}
		& \begin{cases}
			\hat{u}_h = g_D, & \\
			\widehat{\ri\kappa \sigma}_h = - \nabla_h u_h - \ri \kappa {\pa} (u_h-g_D) \bn_K, &  
		\end{cases} & \text{(i)} \\
		& \begin{cases}
			\hat{u}_h = \dlgraffa u_h\drgraffa + \frac12 \bn_K \cdot \bn_\Gamma g_D  - \pb\ri \xi^{-1} \left( \llbracket \nabla_h u_h  \rrbracket_N - g_N \right), & \\
			\widehat{\ri\kappa\sigma}_h = - \dlgraffa \nabla_h u_h\drgraffa - \frac12 \bn_K \cdot \bn_\Gamma \nabla_h g_D - \ri \xi \pa \left( \llbracket u_h  \rrbracket_N -  g_D \bn_\Gamma \right). &  
		\end{cases} & \text{(ii)}
	\end{align*}
	For problem (ii) only, $\xi$ is defined, following \cite{Howarth2014}, on a face $F$ lying on $\Gamma$ as
	\begin{equation} \label{xi_2}
		\xi = {\frac{\Re (\kappa_{\text{o}}) + \Re (\kappa_{\text{i}})}2.}
	\end{equation}
	
	For the TDG scheme in the Dirichlet case, we refer to \cite[Sect. 3]{Kapita2018}. In the transmission problem, substituting the numerical fluxes in (\ref{TDG}) and summing over all the mesh elements, we get the following TDG scheme: Find $u_h \in V_p(\mathcal{M}_h)$ such that for all $v_h \in V_p(\mathcal{M}_h)$
	\begin{equation} \label{bilin_trans}
		\mathcal{A}_h^M(u_h, v_h) =\ell_h(v_h),
	\end{equation}
	where
	\begin{align} \label{DtN_bilinear_exact_trans}
		& \mathcal{A}^M_h(u, v) :=\\ \nonumber
		& \int_{F_h^I} \left( \dlgraffa u\drgraffa \llbracket \overline{\nabla_h v} \rrbracket_N  -  \dlgraffa\nabla_h u\drgraffa \cdot \llbracket \overline{v} \rrbracket_N - \ri \kappa \hspace{0.05cm}  {\pa} \hspace{0.05cm} \llbracket  u\rrbracket_N \cdot \llbracket \overline{v} \rrbracket_N    - \ri \kappa^{-1} {\pb} \hspace{0.05cm} \llbracket \nabla_h u\rrbracket_N  \hspace{0.05cm} \llbracket \overline{\nabla_h v} \rrbracket_N \right) \rd s \\ \nonumber
		& + \int_{\Gamma} \left( \dlgraffa u\drgraffa \llbracket \overline{\nabla_h v} \rrbracket_N  -  \dlgraffa\nabla_h u\drgraffa \cdot \llbracket \overline{v} \rrbracket_N - \ri \xi \hspace{0.05cm}  {\pa} \hspace{0.05cm} \llbracket  u\rrbracket_N \cdot \llbracket \overline{v} \rrbracket_N    -\ri \xi^{-1} {\pb} \hspace{0.05cm} \llbracket \nabla_h u\rrbracket_N  \hspace{0.05cm} \llbracket \overline{\nabla_h v} \rrbracket_N \right) \rd s \\ \nonumber
		& + \int_{\Gamma_{R}} \biggl(  u \hspace{0.05cm} \overline{\partial_{\bn_R} v}  - T_M u \hspace{0.05cm} \overline{v} -{\pd} \hspace{0.05cm} \ri \kappa^{-1} \left(\partial_{\bn_R} u -T_M u\right) \hspace{0.05cm} \overline{\left(\partial_{\bn_R} v -T_M v\right)} \biggl) \rd s,
	\end{align}
	and
	\begin{align} \label{operator_trans}
		\ell_h(v) := &  \int_{\Gamma} \bigg[ - \left( \frac12\bn_K \cdot \bn_\Gamma g_D + \pb \ri \xi^{-1} g_N \right) \llbracket \overline{\nabla_h v} \rrbracket_N \\ \nonumber
		& + \left( \frac12 \bn_K \cdot \bn_\Gamma \nabla_h g_D - \pa \ri \xi g_D \bn_\Gamma \right) \cdot \llbracket \overline{v} \rrbracket_N \bigg] \rd s.
	\end{align}
	
	The consistency of the numerical fluxes gives the consistency of the method \cite{Arnold2002}.
	Moreover, the quasi-optimality of the numerical method is proven in \cite[Prop. 3]{Kapita2018} for the Dirichlet problem. Similar arguments lead to the same results also for the transmission problem.
	
	In our code, we make use of plane waves as basis function of the TDG method. For a mesh element $K \in \Mh$, we denote by $V_p(K)$ the plane wave space on $K$ spanned by $p$ plane waves, $p \in \mathbb{N}$:
	\begin{equation} \label{pwspace}
		V_p(K) = \bigg\{ \; v \in L^2(K) :  v(\bx) = \sum_{j=1}^p \eta_j 
		\exp\{\ri\kappa\bd_j \cdot \bx \}, \hspace{0.3cm} \eta_j \in \mathbb{C} \; \bigg\},
	\end{equation}
	where $\{ \bd_j \}_{j=1}^p \subset \mathbb{R}^2$, with $|\bd_j|=1$, are different propagation directions.
	To obtain isotropic approximations, uniformly-spaced directions	can be chosen as 
	$\bd_j = ( \cos\frac{2\pi j}p, \sin\frac{2\pi j}p)$, $j=1,\ldots,p$. 
	For simplicity, we choose the same number $p$ of directions in every element $K\in\Mh$.
	The value of $\kappa$ depends on the region where the element $K$ is located;
	recall that we consider meshes such that $n$ and $\kappa$ are constant inside each element. We define the global discrete space $V_p(\Mh)$ as
	\begin{equation} \label{pwglobal}
		V_p(\Mh)=\bigoplus_{K\in\Mh}V_p(K)
		=\big\{ \; v \in L^2(\Omega) : v_{|K} \in V_p(K), \;\forall K \in \Mh \; \big\}.
	\end{equation}
	An important advantage of the use of plane waves is that matrix entries on linear edges can easily be computed analytically, reducing the errors caused by numerical quadrature and the computational effort. We refer to \cite[Sect. 4.3]{Monforte2025} for a detailed description of the derivation of the matrix entries. The integrals over the circular boundary $\GR$ can't be computed analytically, so we rely on a Gauss quadrature to compute them; we refer to \cite[Sect. 5]{Kapita2018} for a description of this part.
	
	{
		\begin{rem}[Plane waves approximation properties]
			We refer to \cite[Sect. 4]{Kapita2018} for plane waves estimates for the DtN-TDG method.
			Denoting with $u^M$ the solution of the truncated version of \eqref{bvp}, that is the one we approximate numerically, the error writes as
			\begin{equation*}
				\| u - u_h \|_{L^2(\Omega)} \leq \| u - u^M \|_{L^2(\Omega)} + \| u^M - u_h \|_{L^2(\Omega)}.
			\end{equation*}
			The first term in the inequality is the truncation error of the truncated DtN operator, and it has been studied in \cite{Koyama2009}; see \cite[Prop. 4.1]{Koyama2009} for an explicit estimate of $\| u - u^M \|_{L^2(\Omega)}$.
			
			The $\| u^M - u_h \|_{L^2(\Omega)}$ term is the error of the DtN-TDG method. It has been shown \cite[Thm.~5.2, Cor.~5.5]{Moiola2011} that plane wave spaces admit best-approximation bounds that converge exponentially in $p$ and algebraically in $h$. We refer to \cite[Thm. 5]{Kapita2018} for explicit plane wave estimates for the DtN-TDG method.
		\end{rem}
	}
	
	\section{TMATDG: combining TDG with TMATROM}
	Using the \textsc{tmatrom} package \cite{Ganesh2017}, we are able to approximate the T-matrix of a polygonal scatterer with the DtN-TDG method. Given a wavenumber $\kappa$ and a scatterer $D$, which is described by its geometry and could either be penetrable or impenetrable by waves, we use the DtN-TDG method to solve multiple scattering problems, using the regular wavefunctions as incident fields. We then use the solution obtained to get the far-field and compute the approximation of the T-matrix $\bT$. Once computed, $\bT$ allows to easily solve the scattering problem for any incident field $u^\inc$ and to deal with translations and rotations of the scatterer. To compute the far-field of $u^\scat$ we can use the formula \cite{ColtonKress2019}:
	\begin{equation} \label{eq:farfieldcomp}
		u_\infty^\scat(\theta) = \frac{\re^{\ri\frac\pi4}}{\sqrt{8\pi\kappa}}\int_{\Sigma} \left( u(\bx) \dn \re^{-\ri \kappa \bx \cdot \bd} - \dn u(\bx) \dn \re^{-\ri \kappa \bx \cdot \bd} \right) \rd s(\bx),
	\end{equation}
	where $\bd = (\cos \theta, \sin \theta)$ and $\Sigma$ is a closed curve on which $u^\scat$ admits Dirichlet and Neumann traces. In our setting it is convenient to choose $\Sigma = \GR$.
	
	A core part of the \textsc{tmatdg} package is the approximation of multiple scattering problems. When dealing with an ensemble of scatterers, we have to distinguish between the shape and the orientation, i.e. the position and rotation angle, of a scatterer. An ensemble of obstacles can be composed of multiple copies of the same scatterer that differ only by the orientation: in this case it is convenient to compute the T-matrix of that scatterer only once and then deal with the change of orientation using the addition theorem and the T-matrix properties. This reduces drastically the computational effort. 
	
	\subsection{Approximating the T-matrix of a single obstacle}
	As we said, the \textsc{tmatrom} package allows to include a newly defined numerical method in the approximation of the far-field of a given regular wavefunction; in the \textsc{tmatdg} package we include the DtN-TDG method using the \texttt{TDGsolver} class. We refer to the \textsc{tmatrom} manual for a complete description of how a solver is used to approximate a T-matrix.
	
	Since we are dealing with polygonal scatterers, we identify the scatterer shape by its vertices; moreover, since we deal both with penetrable and impenetrable scatterers, also the type of the scatterer and its refraction index $\nI$ are a parameter. Other parameters are the wavenumber $\kappa$ and the DtN-TDG parameters, i.e. the mesh width $h$ and the number of plane wave directions $p$. The number $M$ of Hankel functions in the truncation of the DtN operator \eqref{DtN} is chosen --- according to \cite[Sect. 6]{Kapita2018} --- so that $M > \kappa R$, where $R$ is the radius of the circular boundary $\GR$. 
	
	We now briefly describe the main functions and classes for solving a single scattering problem:
	\begin{itemize}
		\item The \texttt{TDGsolver} class takes as input the problem parameters, builds the linear system for the DtN-TDG method and solves it for any given incident regular wavefunction, computing the far-field using \eqref{eq:farfieldcomp};
		\item The \texttt{ComputeTMatrix} function takes as input the problem parameters and returns a T-matrix and an instance of the \texttt{TDGsolver} class corresponding to the scatterer; the expansion order (i.e.\ the T-matrix dimension) is chosen accordingly to the theory. For simplicity, the T-matrix is always computed with the scatterer centered in the origin and then it is translated to the scatterer center, which is computed from the geometry;
		\item The \texttt{PlotSolution} function solves the scattering problem for a given incident field $u^\inc$, using the T-matrix approximation, and plots the scattered or total field in a given region. We recall that the expansion \eqref{eq:expansion} holds only outside of the scatterer radius; to plot the field near the scatterer, we use the DtN-TDG method to compute the solution;
		\item The \texttt{RotateTmat} function computes the T-matrix and the solver corresponding to the scatterer rotated by a given rotation angle; the T-matrix is rotated according to \eqref{eq:rotation}, while the domain mesh is simply rotated.
	\end{itemize}
	
	\subsection{Dealing with multiple scattering problems}
	We focus on configurations in which $N$ obstacles assume only $N_S \leq N$ shapes; every obstacle carries its shape, position and rotation angle. The computational cost of approximating the T-matrix for every shape is small and the process is easily done in parallel. After computing the $N_S$ T-matrices, they are translated and rotated to solve the multiple scattering problem.
	
	We describe the functions used to solve a multiple scattering problem:
	\begin{itemize}
		\item The \texttt{MultiScatt} function takes in input a cell array containing the information on the geometry and the type of every shape, and a list with three $N$-dimensional arrays describing the obstacles arrangement: the shape (represented by a vector of integers between $1$ and $N_S$), the position (a complex-valued vector), and  the rotation angle. \texttt{MultiScatt} computes the T-matrices for every shape, and then uses the \textsc{tmatrom} package to solve a multiple scattering problem iteratively with GMRES, finally plotting the total field. It has also the possibility to plot the solution near the scatterers: this is done solving a scattering problem with the DtN-TDG method on the circle using the T-matrix solution as incident field; this can be done only if the circles $\GR$ enclosing the obstacles are not intersecting. The function has the option to save the matrices and solvers created for later uses;
		\item The \texttt{MultiTmatSolve} function allows to solve a multiple scattering when the T-matrices have already been computed, without the need to approximate them again; this reduces the computational cost and time, allowing to solve different scattering problems with a new arrangement or incident field.
	\end{itemize}
	
	\subsection{Examples}
	
	\subsubsection{Single scatterer, rotation and translation} \label{ex:single}
	We show how to approximate the T-matrix of a single polygonal scatterer $D$ with $V$ vertices, and how to use and modify it. An obstacle is represented by a structure with 3 fields: the polygon vertices (a $V \times 2$ matrix), the scatterer type (either \texttt{'dir'} or \texttt{'trans'}), and the value of $\nI$ (for penetrable scatterers). We also fix the wavenumber and the TDG parameters.
	To compute the approximate T-matrix we define the parameters and call the \texttt{ComputeTMatrix} function:
	\begin{lstlisting}[language=Octave]
		k = 5; h = 0.5; p = 20; 
		scatt.vertices = [-1, -1; -1, 1; 1, 1; 1, -1]; 
		scatt.type = 'trans'; 
		scatt.n_in = 3+1i;
		[tmat,solver]=ComputeTMatrix(k,h,p,scatt);
	\end{lstlisting}
	Then, we choose an incident plane wave $u^\inc$ of angle $\theta$ and solve and plot the solution on a given domain using the \texttt{PlotSolution} function; we can specify if we want to plot also near the scatterer:
	\begin{lstlisting}[language=Octave]
		theta = -pi/3; uinc = plane_wave(theta,k); 
		PlotPar.type_plot = 'tot'; PlotPar.inside = true; 
		PlotPar.limX = [-5,5]; PlotPar.limY = [-5,5];
		PlotSolution(tmat,solver,uinc,PlotPar);
	\end{lstlisting}
	\begin{figure}[htbp]
		\centering
		\includegraphics[width=.3\textwidth]{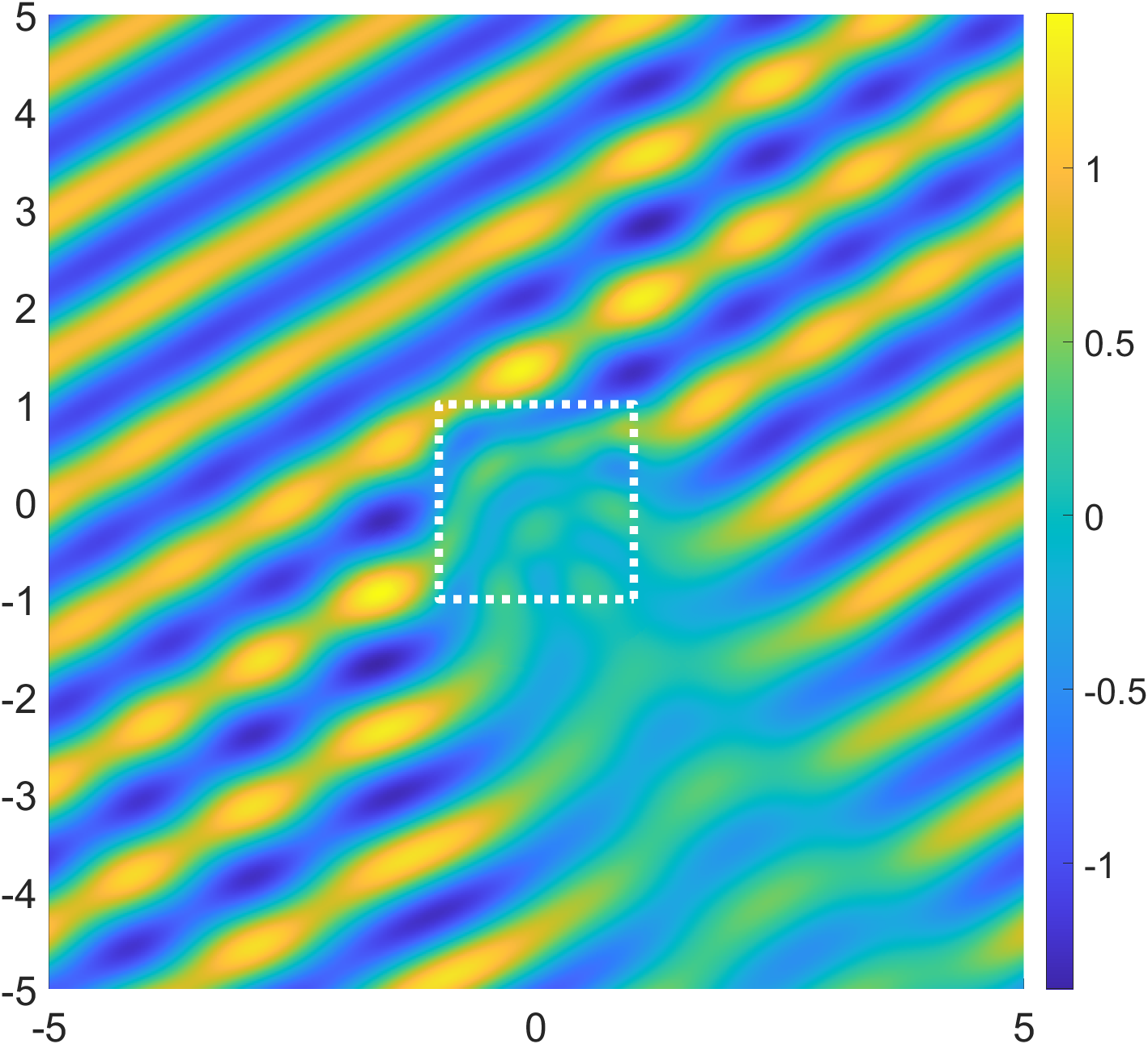}\hfil\includegraphics[width=.3\textwidth]{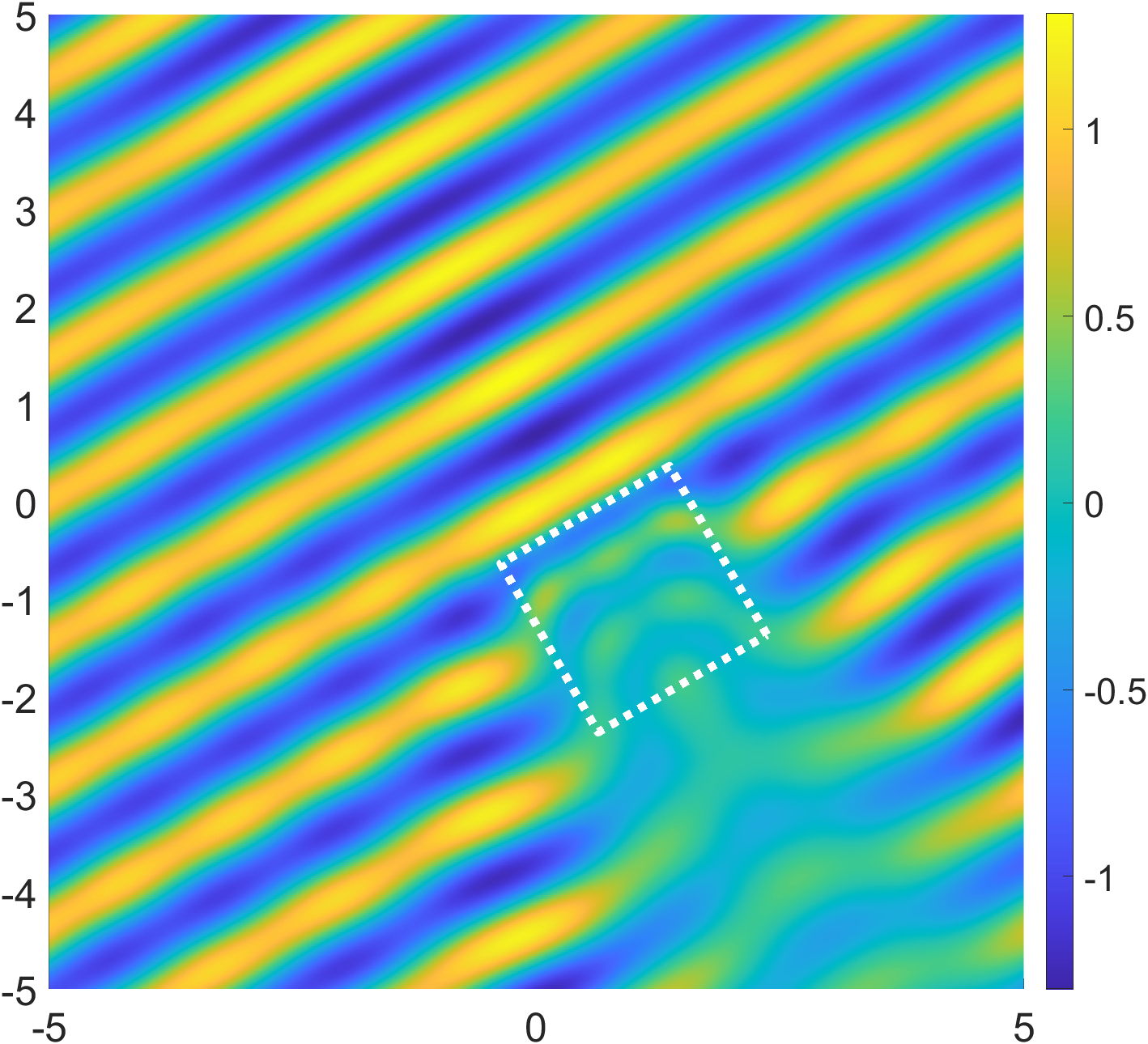}\hfil\includegraphics[width=.3\textwidth]{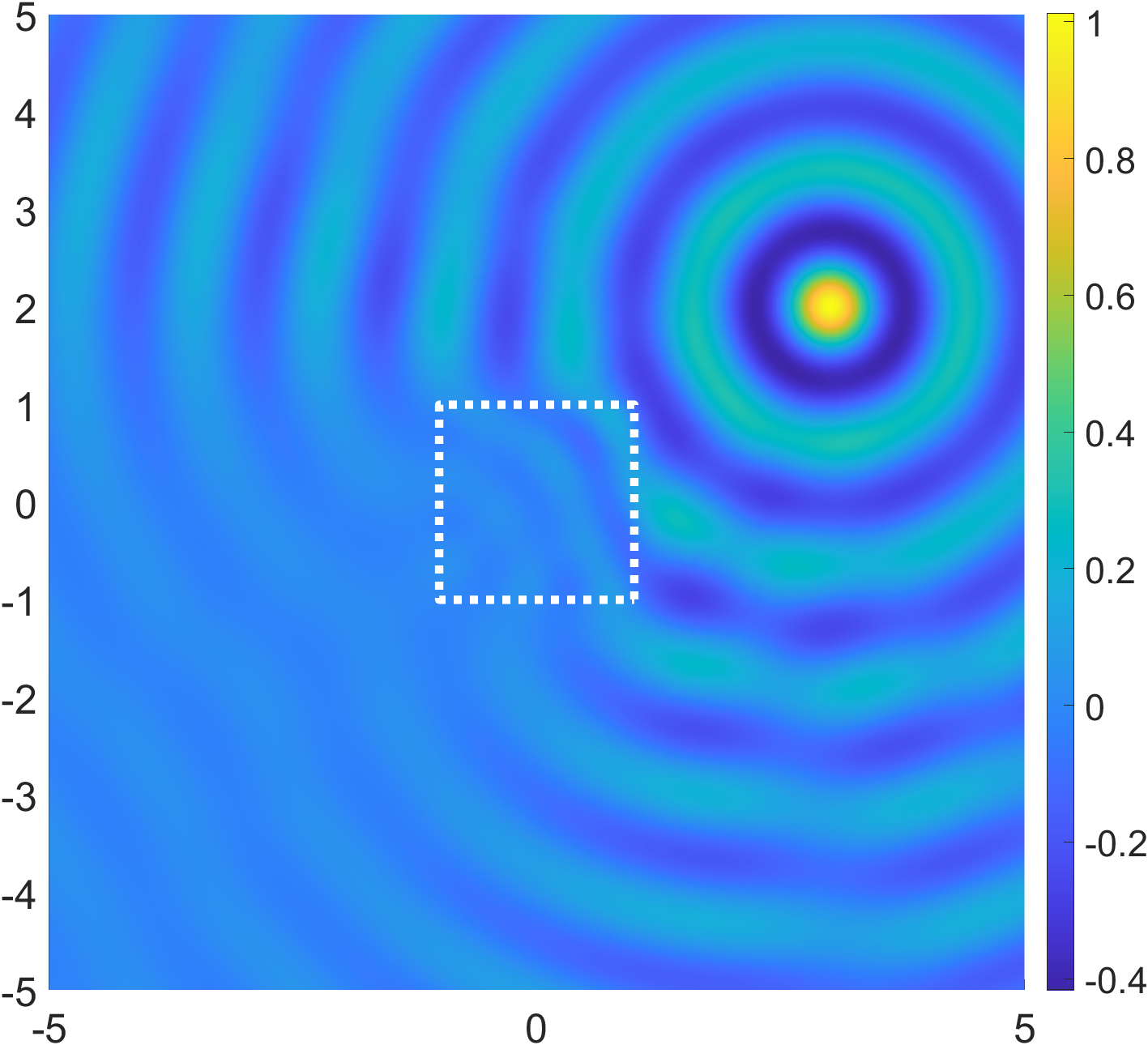}
		\caption{Real part of total field for the square penetrable scatterer of Section \ref{ex:single} in 3 different cases: Left: plane incident wave; Center: rotated and translated scatterer with the same incident wave; Right: circular incident wave.}
	\end{figure}
	We can rotate the scatterer by a given angle using \texttt{RotateTmat} and translate it using the \textsc{tmatrom} \texttt{setOrigin} function; then, we can solve the new scattering problem with the same $u^\inc$:
	\begin{lstlisting}[language=Octave]
		rotation_angle = pi/6;
		[rotTmat,rotSolver] = RotateTmat(tmat,solver,rotation_angle);
		rotTmat.setOrigin(1-1i); 
		PlotSolution(rotTmat,rotSolver,uinc,PlotPar);
	\end{lstlisting}
	It is easy to change $u^\inc$; we choose a circular wave centered in $(3,2)$:
	\begin{lstlisting}[language=Octave]
		uinc_new = point_source(3+2i,k);  
		PlotSolution(tmat,solver,uinc_new,PlotPar);
	\end{lstlisting}
	
	\subsubsection{Multiple scattering problem and configuration changes} \label{ex:multiple}
	We solve the multiple scattering problem with $N_S = 3$ shapes: an impenetrable square, a penetrable triangle with $\nI = 2.5$ and an impenetrable cross-shaped obstacle; the information are stored in the \texttt{ScatShape} cell array:
	\begin{lstlisting}[language=Octave]
		NShape = 3; ScatShape=cell(NShape,1);
		ScatShape{1}.vertices = [1/3, 1/3; 1/3, 1; -1/3, 1; -1/3, 1/3; -1, 1/3;
		-1, -1/3; -1/3, -1/3; -1/3, -1; 1/3, -1; 1/3, -1/3; 1, -1/3; 1, 1/3]; 
		ScatShape{1}.type = 'dir';
		ScatShape{2}.vertices = [0, 1; -sqrt(3)/2, -1/2; sqrt(3)/2, -1/2,]; 
		ScatShape{2}.type = 'trans'; ScatShape{2}.n_in = 2.5;
		ScatShape{3}.vertices = [1, 0; 0, 1; -1, 0; 0,-1]; 
		ScatShape{3}.type = 'dir';
	\end{lstlisting}
	In \texttt{ScatArr} we describe the ensemble arrangement; we have $N=5$ obstacles:
	\begin{lstlisting}[language=Octave]
		ScatArr.shape = [1; 1; 2; 3; 2]; 
		ScatArr.pos = [-4-4i; 4-3.5i; 0; -3+4i; 3.5+3i];
		ScatArr.rot = [-pi/4; 0; 0; 0; pi]; 
	\end{lstlisting}
	A plane wave with wavenumber $\kappa = 10$ and direction $(-1,1)$ is incident on the ensemble; we use the \texttt{MultiScatt} function to solve the multiple scattering problem and plot the solution, storing the T-matrices and solvers for later use:
	\begin{lstlisting}[language=Octave]
		k = 10; h = 0.5; p = 20;
		theta = 3*pi/4; uinc = plane_wave(theta,k); 
		PlotPar.inside = true; PlotPar.limX=[-7,7]; PlotPar.limY=[-7,7];
		SavePath.choice=true; SavePath.file = 'MultiTest.mat';
		MultiScatt(k,h,p,uinc,ScatShape,ScatArr,PlotPar,SavePath)
	\end{lstlisting}
	\begin{figure}[htbp]
		\centering
		\includegraphics[height=.4\textwidth]{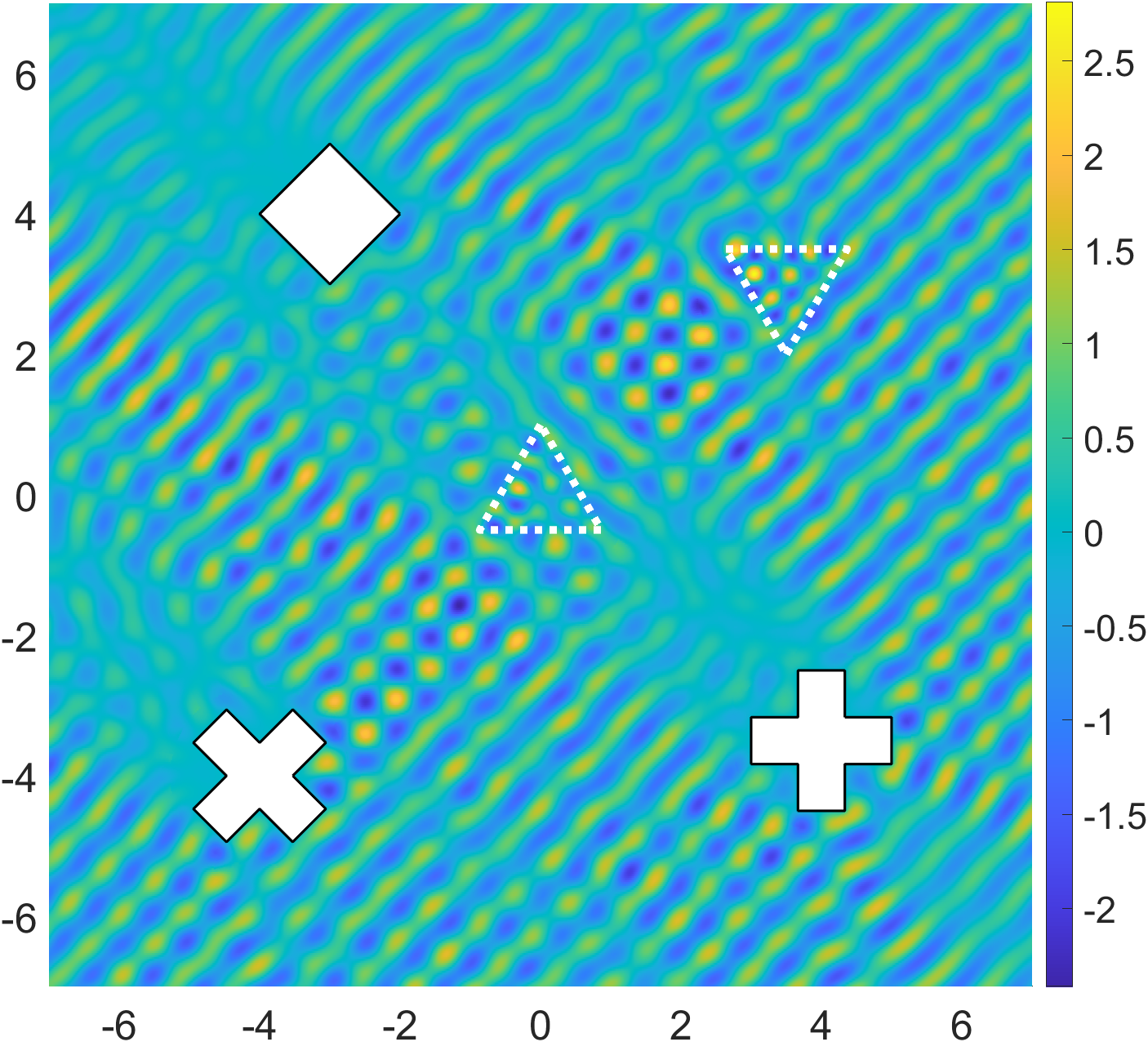}\hfil\includegraphics[height=.4\textwidth]{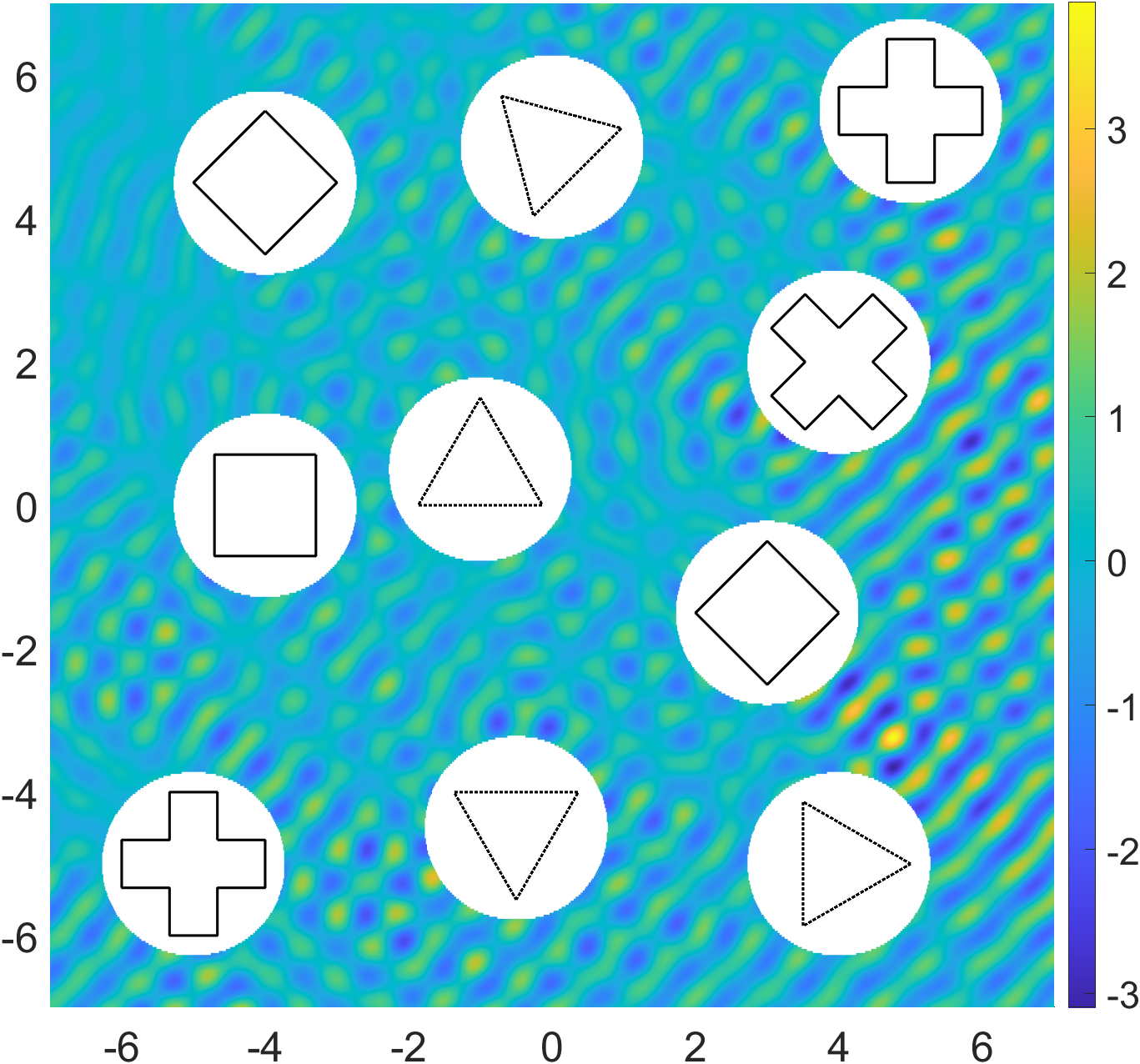}
		\caption{Real part of total field for the ensemble of scatterers of Section \ref{ex:multiple} with two different arrangements. On the left we plot near the obstacles; on the right this is not possible since the obstacles are too close to each other and the TDG solutions are computed on disks of radius $R_D+2h$.}
		\label{fig:multiple}
	\end{figure}
	Then, we change the number of obstacles to $N=10$ and modify the arrangement; we use the \texttt{MultiTmatSolve} function since we already saved the T-matrices and solvers:
	\begin{lstlisting}[language=Octave]
		MultiTmatSolve(k,uinc,tmat,solver,ScatArrNew,PlotPar);
	\end{lstlisting}
	In Figure \ref{fig:multiple} we plot the real part of the total field for both the examples. 
	
	{For the first arrangement, we perform a $p$- and $M$-convergence test, computing the $L^2$-norm of the relative error against a refined solution on the domain $[-7,7]\times[-7,7]$, showing that both the DtN-TDG method and the T-matrix approximation are improving. First, we fix the mesh width $h=0.5$ and the DtN truncation order $M=20$ and increase only the number $p$ of basis functions inside every element of the mesh, computing the error against a refined solution with $p=25$. In the second case, we fix $h=0.5$ and $p=20$ and increase the Fourier order of truncation $M$, computing the error against a refined solution with $M=50$.
		The results are displayed in Figure \ref{fig:conv}: the $p$-convergence is in agreement with other DtN-TDG results in presence of corners, e.g.\ \cite[Fig.~11]{Monforte2025}, while the $M$-convergence is in agreement with the results in \cite[Fig.~2]{Kapita2018}.
		\begin{figure}[htbp]
			\centering
			\includegraphics[height=.3\textwidth]{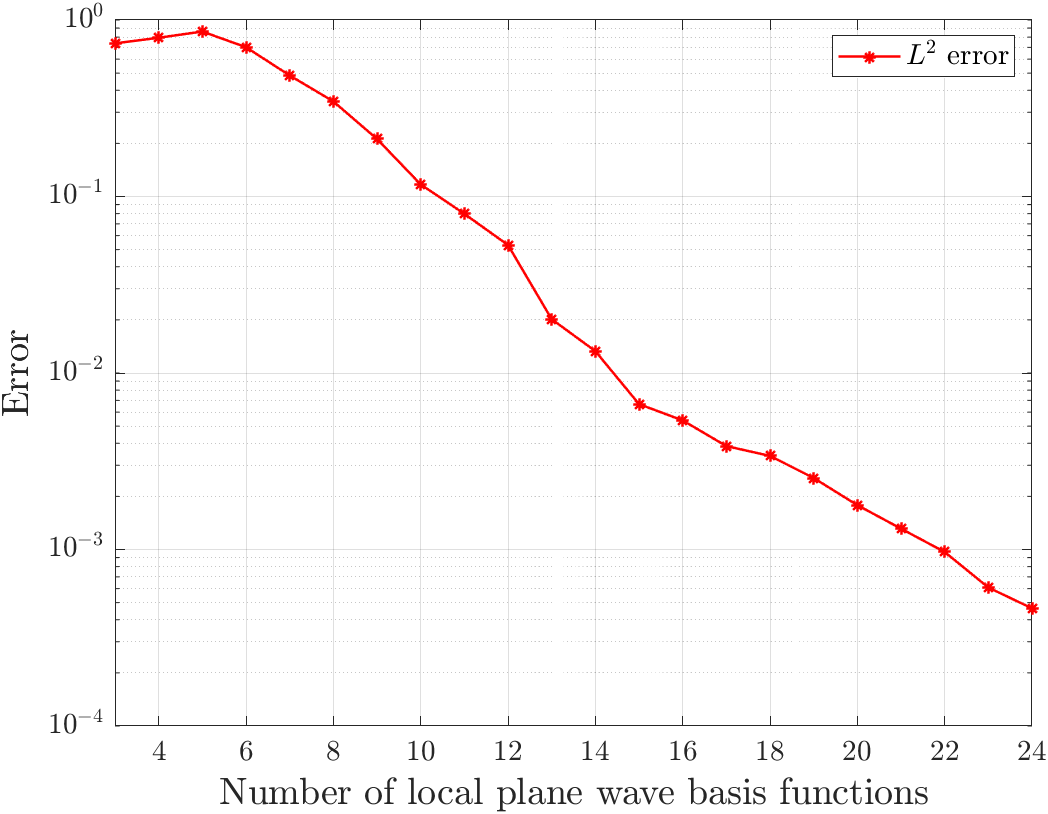}\hfil\includegraphics[height=.3\textwidth]{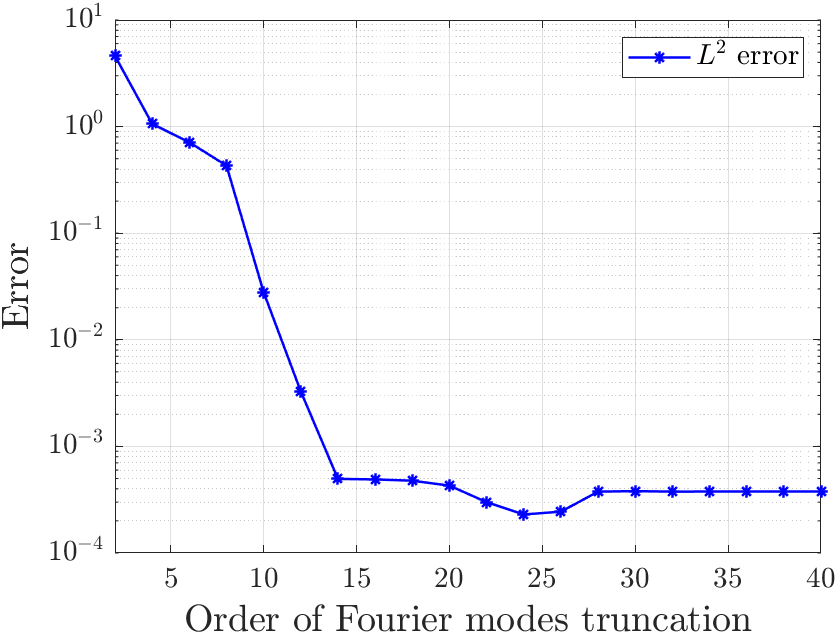}
			\caption{Convergence tests for the ensemble of scatterers of Section \ref{ex:multiple}. Left: the $p$-convergence of the relative $L^2$ error, computed on $[-7,7]\times[-7,7]$, for $p \in \{3,\ldots,24\}$. Right: $M$-convergence of the relative $L^2$ error, computed on the same domain, for $M \in \{2,\ldots,40\}$.}
			\label{fig:conv}
		\end{figure}
	}
	
	\subsubsection{A parameter-depending problem} \label{ex:parameter}
\begin{figure}[htbp]
		\centering
		\includegraphics[height=.3\textwidth]{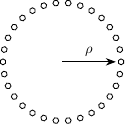}\hfil\includegraphics[height=.3\textwidth]{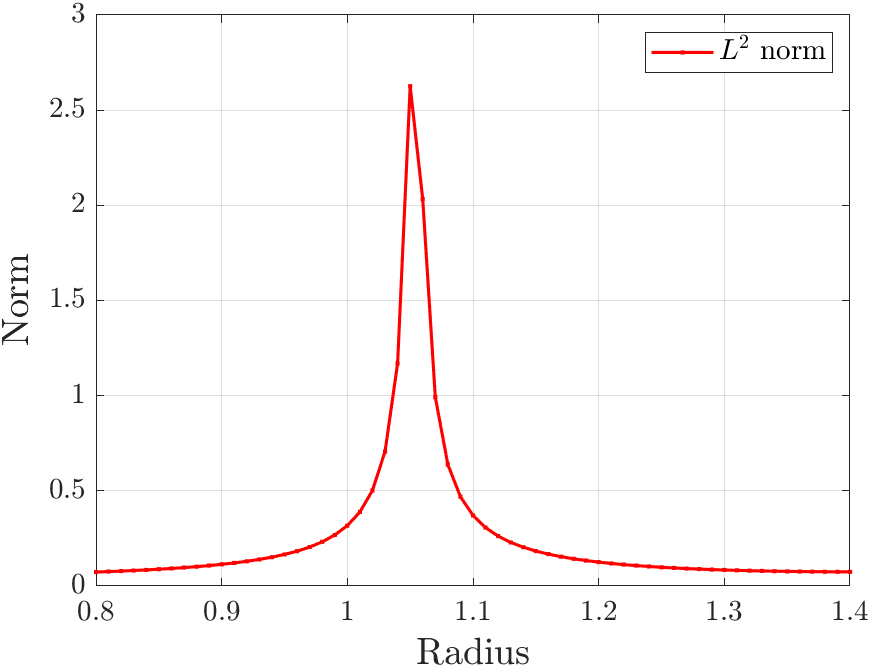}
		\caption{Left: the obstacle arrangement of Section \ref{ex:parameter}. Right: the $L^2$-norm of the solution computed on the ball of radius $0.5$ for different values of $\rho$.}
		\label{fig:disp_norm}
	\end{figure}

	We arrange $N=30$ copies of an hexagonal impenetrable obstacle as in Figure \ref{fig:disp_norm}: the hexagons centers are located on a circle of radius $\rho$, which may vary; each hexagon is inscribed in a circle of radius $R_D = 0.05$. Every copy of the obstacle is just a translation and rotation of the same reference shape, such that a side faces the origin. 
	
	We fix a wavenumber $\kappa^*=2.39$ and an incident field $u^\inc$, a circular wave centered in $(2,0)$; the choice of $\kappa^*$ is such that $J_0(\kappa^*) \approx 0$, so we expect to see a resonance when the obstacles are disposed around the unit circle \cite{Hewett2016}.
	We vary the radius $\rho$ of the circle between $0.8$ and $1.4$ and compute the $L^2$-norm of the total field on the ball of radius $0.5$ centered at the origin. We aim to find the value of $\rho$ for which the norm is maximal. In this setting, using the T-matrix method is convenient, since the matrix is computed only once and then reused for all the values of $\rho$, changing only the arrangement. 
	
	The results are displayed in Figure \ref{fig:disp_norm}: we observe that we reach a maximum for $\rho = 1.05$; this is explained by the fact that the internal circle formed by the hexagons has radius $1$, so we get the resonance.
	\begin{figure}[htbp]
		\centering
		\includegraphics[height=.37\textwidth]{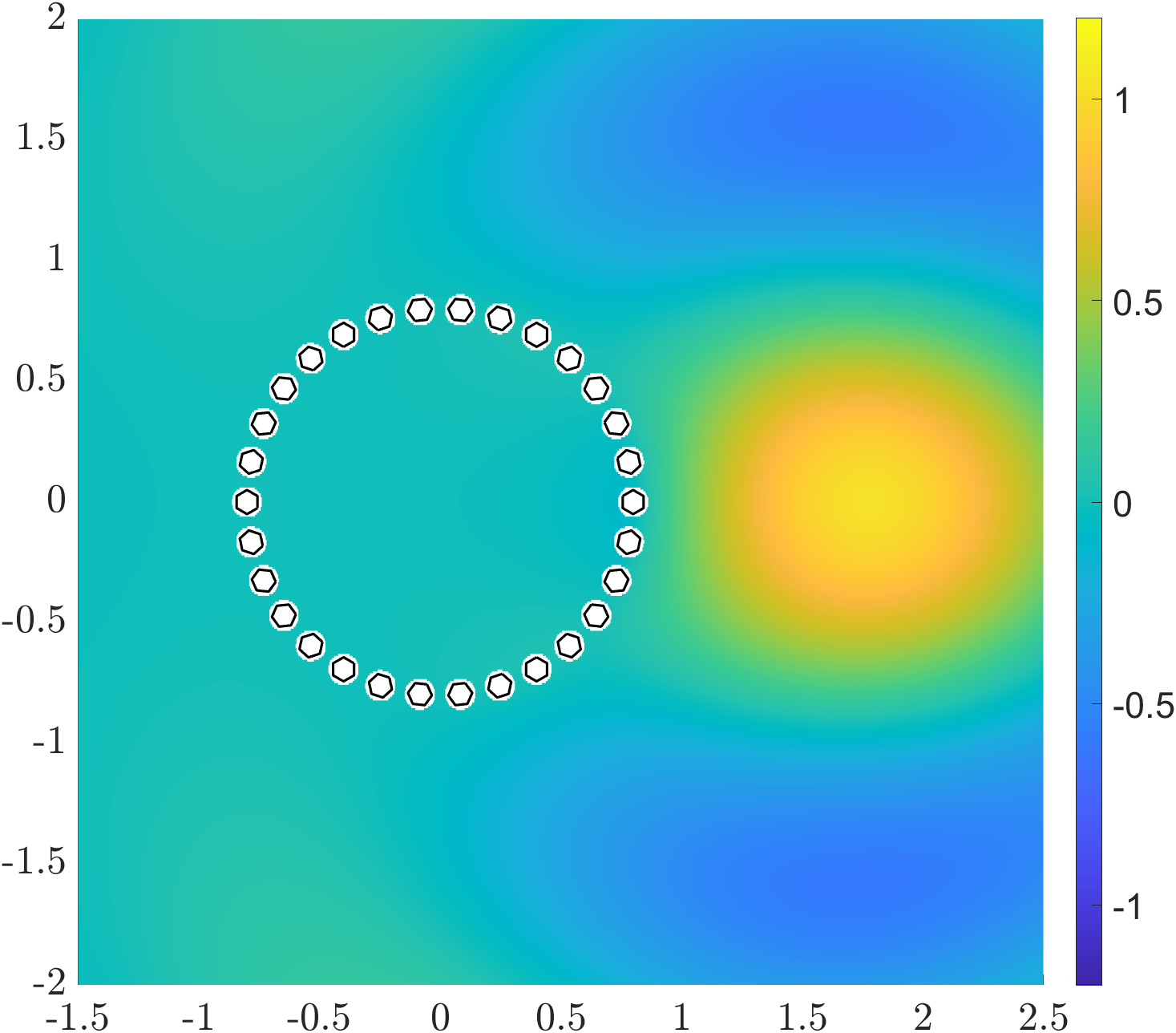}\hfil\includegraphics[height=.37\textwidth]{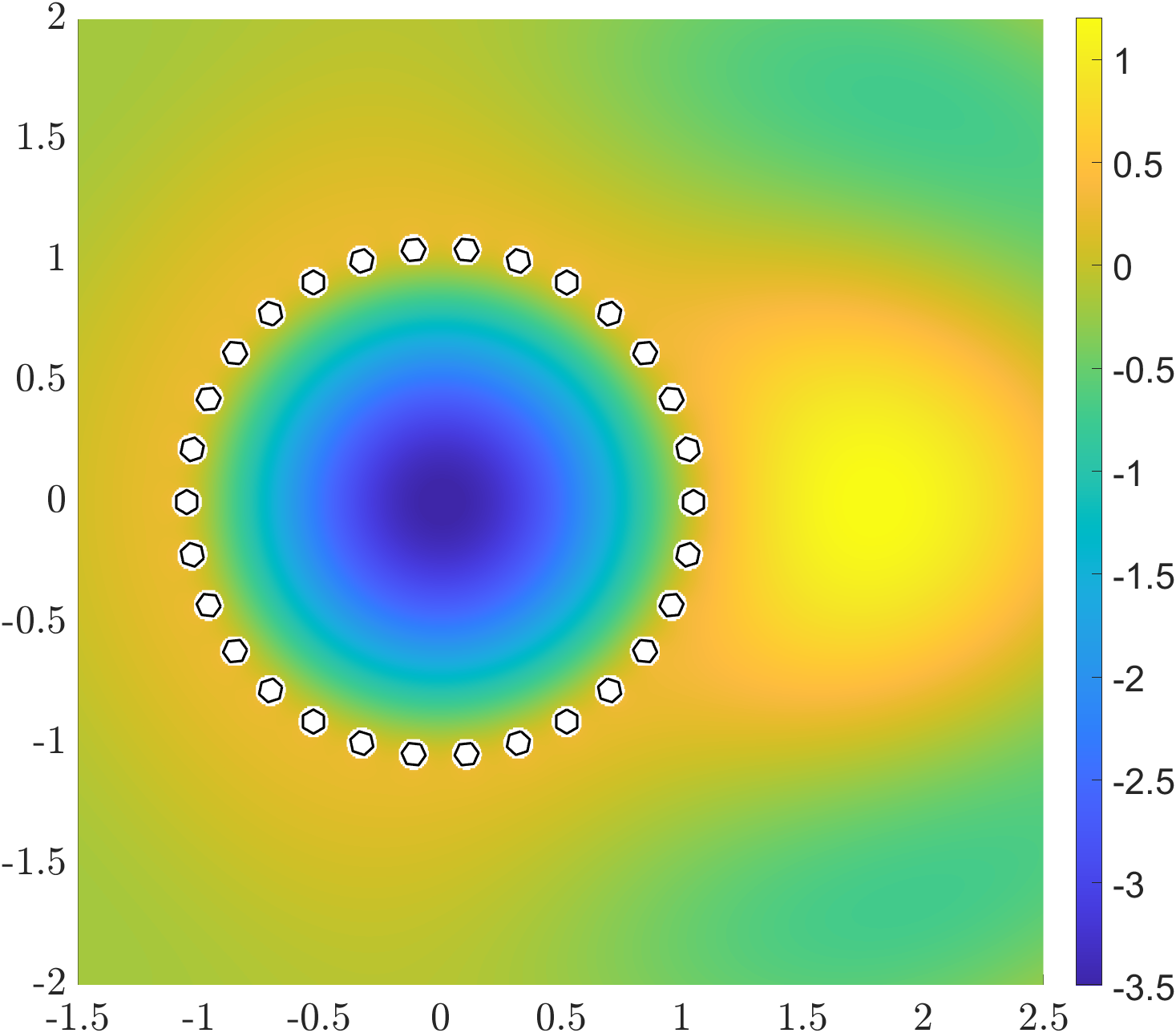}
		\caption{Real part of total field for the scatterer arrangement of Section \ref{ex:parameter}, with $\rho = 0.8$ and $\rho = 1.05$ respectively.}
		\label{fig:plot_hexagons}
	\end{figure}
	
	In Figure \ref{fig:plot_hexagons} we display the real part of the total field for $\rho=0.8$ and $\rho=1.05$: as we can observe, the absolute value of the solution inside the ball is much bigger if we choose the resonant radius, while in the other case there is no resonance.
	
	\section*{Declarations}
	\textbf{Acknowledgements.} The author thanks Andrea Moiola (University of Pavia) for providing valuable tips and guidance during the preparation of this article, and the anonymous referee for providing many constructive comments that improved the clarity of the paper. The author is also member of the GNCS-INdAM group.
	
	\medskip\noindent
	\textbf{Funding.} The author was financially supported by the PRIN project ``ASTICE'' (202292JW3F), funded by the European Union--NextGenerationEU.
	
	\medskip\noindent
	\textbf{Declaration of competing interest.} The author has no relevant financial or non-financial interests to disclose.
	
	\medskip\noindent
	\textbf{Code availability.} The package implemented is available at \url{https://github.com/Arma99dillo/TMATDG}
	
	\printbibliography

@article {Hiptmair2011,
    AUTHOR = {Hiptmair, Ralf and Moiola, Andrea and Perugia, Ilaria},
     TITLE = {Plane wave discontinuous {G}alerkin methods for the 2{D}
              {H}elmholtz equation: analysis of the {$p$}-version},
   JOURNAL = {SIAM J. Numer. Anal.},
  FJOURNAL = {SIAM Journal on Numerical Analysis},
    VOLUME = {49},
      YEAR = {2011},
    NUMBER = {1},
     PAGES = {264--284},
      ISSN = {0036-1429,1095-7170},
   MRCLASS = {65N30 (35J05 65N15 76Q05)},
  MRNUMBER = {2783225}
}

@article {Kapita2018,
    AUTHOR = {Kapita, Shelvean and Monk, Peter},
     TITLE = {A plane wave discontinuous {G}alerkin method with a
              {D}irichlet-to-{N}eumann boundary condition for the scattering
              problem in acoustics},
   JOURNAL = {J. Comput. Appl. Math.},
  FJOURNAL = {Journal of Computational and Applied Mathematics},
    VOLUME = {327},
      YEAR = {2018},
     PAGES = {208--225},
      ISSN = {0377-0427,1879-1778},
   MRCLASS = {65N21 (35J05 35P25 49M25 76Q05)},
  MRNUMBER = {3683156},
MRREVIEWER = {Erwin\ Stein}
}

@phdthesis{Howarth2014,
    author = {Howarth, Charlotta Jasmine},
    title = {New Generation Finite Element Methods For Forward Seismic Modelling},
    school = {University of Reading},
    year = {2014},
    note = {Available on \url{https://www.reading.ac.uk/maths-and-stats/publications/theses-and-dissertations/mathematics-phd-theses}}
}

@Article{Arnold2002,
  author = {D. N. Arnold and F. Brezzi and B. Cockburn and L. D. Marini},
  title = {Unified analysis of discontinuous {G}alerkin methods for elliptic  problems},
  journal = {SIAM J. Numer. Anal.},
  fjournal = {SIAM Journal of Numerical Analysis},
  year = {2002},
  volume = {39},
  number = {5},
  pages = {1749--1779}
}

@Book{ColtonKress2019,
 Author = {Colton, David and Kress, Rainer},
 Title = {Inverse acoustic and electromagnetic scattering theory},
 Edition = {4th expanded edition},
 FSeries = {Applied Mathematical Sciences},
 Series = {Appl. Math. Sci.},
 ISSN = {0066-5452},
 Volume = {93},
 ISBN = {978-3-030-30350-1; 978-3-030-30351-8},
 Year = {2019},
 Publisher = {Cham: Springer},
 Language = {English},
 DOI = {10.1007/978-3-030-30351-8},
 zbMATH = {7114984},
 Zbl = {1425.35001}
}

@ARTICLE{Waterman1965,
	author={Waterman, P.C.},
	journal={Proceedings of the IEEE}, 
	title={Matrix formulation of electromagnetic scattering}, 
	year={1965},
	volume={53},
	number={8},
	pages={805-812},
	doi={10.1109/PROC.1965.4058}}

@article{Waterman1971,
	title = {Symmetry, Unitarity, and Geometry in Electromagnetic Scattering},
	author = {Waterman, P. C.},
	journal = {Phys. Rev. D},
	volume = {3},
	issue = {4},
	pages = {825--839},
	numpages = {0},
	year = {1971},
	month = {Feb},
	publisher = {American Physical Society},
	doi = {10.1103/PhysRevD.3.825},
	url = {https://link.aps.org/doi/10.1103/PhysRevD.3.825}
}

@article{Ganesh2009,
	title = "A far-field based T-matrix method for two dimensional obstacle scattering",
	author = "M. Ganesh and Hawkins, {S. C.}",
	year = "2009",
	doi = "10.21914/anziamj.v51i0.2581",
	language = "English",
	volume = "51",
	pages = "C215--C230",
	journal = "ANZIAM Journal",
	issn = "1446-1811",
	publisher = "Cambridge University Press (CUP)",
}

@article{Wiscombe1980,
	author = {W. J. Wiscombe},
	journal = {Appl. Opt.},
	number = {9},
	pages = {1505--1509},
	publisher = {Optica Publishing Group},
	title = {Improved Mie scattering algorithms},
	volume = {19},
	month = {May},
	year = {1980},
	url = {https://opg.optica.org/ao/abstract.cfm?URI=ao-19-9-1505},
	doi = {10.1364/AO.19.001505},
}

@article{Ganesh2012,
	author = {Ganesh, Mogulapally and Hawkins, S.C. and Hiptmair, Ralf},
	year = {2012},
	month = {10},
	pages = {},
	title = {Convergence analysis with parameter estimates for a reduced basis acoustic scattering T-matrix method},
	volume = {32},
	journal = {IMA Journal of Numerical Analysis},
	doi = {10.1093/imanum/drr041}
}

@article{Waterman1969,
	author = {Waterman, P. C.},
	title = {New Formulation of Acoustic Scattering},
	journal = {The Journal of the Acoustical Society of America},
	volume = {45},
	number = {6},
	pages = {1417-1429},
	year = {1969},
	month = {06},
	issn = {0001-4966},
	doi = {10.1121/1.1911619},
	url = {https://doi.org/10.1121/1.1911619},
	eprint = {https://pubs.aip.org/asa/jasa/article-pdf/45/6/1417/18764237/1417_1_online.pdf},
}

@article{Mishchenko1994,
	title = {T-matrix computations of light scattering by large spheroidal particles},
	journal = {Optics Communications},
	volume = {109},
	number = {1},
	pages = {16-21},
	year = {1994},
	issn = {0030-4018},
	doi = {https://doi.org/10.1016/0030-4018(94)90731-5},
	url = {https://www.sciencedirect.com/science/article/pii/0030401894907315},
	author = {Michael I. Mishchenko and Larry D. Travis}
}

@article{Mishchenko1996,
	title = {T-matrix computations of light scattering by nonspherical particles: A review},
	journal = {Journal of Quantitative Spectroscopy and Radiative Transfer},
	volume = {55},
	number = {5},
	pages = {535-575},
	year = {1996},
	note = {Light Scattering by Non-Spherical Particles},
	issn = {0022-4073},
	doi = {https://doi.org/10.1016/0022-4073(96)00002-7},
	url = {https://www.sciencedirect.com/science/article/pii/0022407396000027},
	author = {Michael I. Mishchenko and Larry D. Travis and Daniel W. Mackowski}
}

@article{Farafonov2010,
	title = {Near- and far-field light scattering by nonspherical particles: Applicability of methods that involve a spherical basis},
	journal = {Optics and Spectroscopy},
	volume = {109},
	number = {3},
	pages = {432-443},
	year = {2010},
	doi = {https://doi.org/10.1134/S0030400X10090195},
	author = {Farafonov, V. G. and Il’in, V. B. and Vinokurov, A. A.}
}

@article{HellmersJens2011,
	author = {Hellmers, Jens and Schmidt, Vladimir and Wriedt, Thomas},
	copyright = {2011 Elsevier Ltd},
	issn = {0022-4073},
	journal = {Journal of quantitative spectroscopy and radiative transfer},
	language = {eng},
	number = {11},
	pages = {1679-1686},
	publisher = {Elsevier Ltd},
	title = {Improving the numerical stability of T-matrix light scattering calculations for extreme particle shapes using the nullfield method with discrete sources},
	volume = {112},
	year = {2011},
}

@inProceedings{Ganesh2008,
	author = "M. Ganesh
	and S. C. Hawkins",
	title = "A far-field based T-matrix method for three dimensional acoustic scattering",
	series = "ANZIAM J.",
	volume = "50",
	pages = "C121--C136",
	year = 2008,
	booktitle = " Proceedings of the 14th Biennial Computational Techniques and Applications Conference, CTAC-2008",
	editor = "Geoffry N. Mercer  and A. J. Roberts",
	month = oct,
	url = http://anziamj.austms.org.au/ojs/index.php/ANZIAMJ/article/view/1441,
	subjclass = "",
}

@article{Ganesh2017,
	author = {Ganesh, M. and Hawkins, S. C.},
	title = {Algorithm 975: TMATROM—A T-Matrix Reduced Order Model Software},
	year = {2017},
	issue_date = {March 2018},
	publisher = {Association for Computing Machinery},
	address = {New York, NY, USA},
	volume = {44},
	number = {1},
	issn = {0098-3500},
	url = {https://doi.org/10.1145/3054945},
	doi = {10.1145/3054945},
	journal = {ACM Trans. Math. Softw.},
	month = jul,
	articleno = {9},
	numpages = {18}
}

@book{Martin2006,
	author = {Martin, P.A},
	title = {Multiple Scattering. Interaction of Time-Harmonic Waves with N Obstacles},
	journal = {Encyclopedia of Mathematics and Its Applications},
	publisher = {Cambridge University Press},
	volume = {107},
	year = {2006},
	isbn = {0521 865549},
}

@article{Ganesh2013,
	title = {A stochastic pseudospectral and T-matrix algorithm for acoustic scattering by a class of multiple particle configurations},
	journal = {Journal of Quantitative Spectroscopy and Radiative Transfer},
	volume = {123},
	pages = {41-52},
	year = {2013},
	note = {Peter C. Waterman and his scientific legacy},
	issn = {0022-4073},
	doi = {https://doi.org/10.1016/j.jqsrt.2013.01.011},
	url = {https://www.sciencedirect.com/science/article/pii/S0022407313000290},
	author = {M. Ganesh and S.C. Hawkins}
}

@article{Monforte2025,
	title={Trefftz Discontinuous Galerkin methods for scattering by periodic structures}, 
	author={Armando Maria Monforte and Andrea Moiola},
	year={2025},
	eprint={2505.23216},
	journal={arXiv: 2505.23216},
	archivePrefix={arXiv},
	primaryClass={math.NA},
	url={https://arxiv.org/abs/2505.23216}, 
}

@article{Hewett2016,
	author = {Hewett, D. P. and Hewitt, I. J.},
	title = {Homogenized boundary conditions and resonance effects in Faraday cages},
	journal = {Proceedings of the Royal Society A: Mathematical, Physical and Engineering Sciences},
	volume = {472},
	number = {2189},
	pages = {20160062},
	year = {2016},
	month = {05},
	issn = {1364-5021},
	doi = {10.1098/rspa.2016.0062},
	url = {https://doi.org/10.1098/rspa.2016.0062},
	eprint = {https://royalsocietypublishing.org/rspa/article-pdf/doi/10.1098/rspa.2016.0062/361313/rspa.2016.0062.pdf},
}

@article{Hawkins2024,
	author = {Hawkins, Stuart C. and Bennetts, Luke G. and Nethercote, Matthew A. and Peter, Malte A. and Peterseim, Daniel and Putley, Henry J. and Verfürth, Barbara},
	title = {Metamaterial applications of Tmatsolver, an easy-to-use software for simulating multiple wave scattering in two dimensions},
	journal = {Proceedings of the Royal Society A: Mathematical, Physical and Engineering Sciences},
	volume = {480},
	number = {2292},
	pages = {20230934},
	year = {2024},
	month = {06},
	issn = {1364-5021},
	doi = {10.1098/rspa.2023.0934},
	url = {https://doi.org/10.1098/rspa.2023.0934},
	eprint = {https://royalsocietypublishing.org/rspa/article-pdf/doi/10.1098/rspa.2023.0934/1229778/rspa.2023.0934.pdf},
}

@article{Koyama2009,
	title = {Error estimates of the finite element method for the exterior Helmholtz problem with a modified DtN boundary condition},
	journal = {Journal of Computational and Applied Mathematics},
	volume = {232},
	number = {1},
	pages = {109-121},
	year = {2009},
	note = {Special Issue: Honor of Professor Hideo Kawarada on the Occasion of his 70th Birthday},
	issn = {0377-0427},
	doi = {https://doi.org/10.1016/j.cam.2008.10.034},
	url = {https://www.sciencedirect.com/science/article/pii/S0377042708005505},
	author = {Daisuke Koyama}
}

@article {Moiola2011,
	AUTHOR = {Moiola, Andrea and Hiptmair, Ralf and Perugia, Ilaria},
	TITLE = {Plane wave approximation of homogeneous {H}elmholtz solutions},
	JOURNAL = {Z. Angew. Math. Phys.},
	FJOURNAL = {Zeitschrift f\"{u}r Angewandte Mathematik und Physik. ZAMP.
	Journal of Applied Mathematics and Physics. Journal de
	Math\'{e}matiques et de Physique Appliqu\'{e}es},
	VOLUME = {62},
	YEAR = {2011},
	NUMBER = {5},
	PAGES = {809--837},
	ISSN = {0044-2275,1420-9039},
	MRCLASS = {35J05 (31B35 35A35 41A30 65N99)},
	MRNUMBER = {2843918},
	MRREVIEWER = {Michele\ Campiti}
}
	
\end{document}